\documentclass[11pt,twoside]{article}
\usepackage{mathrsfs}
\usepackage{amsmath}
\usepackage{amssymb}
\usepackage{fancyhdr}
\usepackage{latexsym}
\usepackage{bbding}
\usepackage{mathrsfs}
\usepackage{wasysym}
\usepackage{cite}
\usepackage{multicol,graphics}

\newtheorem{Th}{Theorem}[section]
\newtheorem{Le}{Lemma}[section]
\newtheorem{De}{Definition}[section]
\newtheorem{Rem}{Remark}[section]
\newtheorem{Cor}{Corollary}[section]

\newcommand{\bth}{\begin{Th}}
\newcommand{\eeth}{\end{Th}}
\newcommand{\ble}{\begin{Le}}
\newcommand{\eele}{\end{Le}}
\newcommand{\bde}{\begin{De}}
\newcommand{\ede}{\end{De}}
\newcommand{\bre}{\begin{Rem}}
\newcommand{\eere}{\end{Rem}}
\newcommand{\bco}{\begin{Cor}}
\newcommand{\eeco}{\end{Cor}}
\newcommand{\p}{\rm{\bf Proof.}}

 \makeatletter
\allowdisplaybreaks
\usepackage{esint}

\newtheorem{lem}{Lemma}[section]

\numberwithin{equation}{section}

\usepackage{authblk}
\title{Riemann-Hilbert problems for bi-axially symmetric null-solutions to iterated perturbed Dirac equations in $\mathbb{R}^{n}$}
\author[a]{Dian Zuo}
\author[b]{Min Ku}
\author[a]{Fuli He\thanks{\text{Corresponding author. E-mail}: hefuli999@163.com.}}
\affil[a]{School of Mathematics and Statistics, Central South University, Changsha 410083, China}
\affil[b]{Institute for Computing and Information Sciences, University of Radboud, the Netherlands}

\date{}

\begin{document}

\maketitle

\begin{abstract}

\noindent
This work addresses Riemann-Hilbert boundary value problems (RHBVPs) for null solutions to iterated perturbed Dirac operators over bi-axially symmetric domains in $\mathbb{R}^n$
with Clifford-algebra-valued variable coefficients. We first resolve the unperturbed case of poly-monogenic functions, i.e., null solutions to iterated Dirac operators, by constructing explicit solutions via a bi-axially adapted Almansi-type decomposition, decoupling hierarchical structures through recursive integral operators. Then, generalizing to vector wave number-perturbed iterated Dirac operators, we extend the decomposition to manage spectral anisotropy while preserving symmetry constraints, ensuring regularity under Clifford-algebraic parameterizations. As a key application, closed-form solutions to the Schwarz problem are derived, demonstrating unified results across classical and higher-dimensional settings. The interplay of symmetry, decomposition, and perturbation theory establishes a cohesive framework for higher-order boundary value challenges in Clifford analysis.
\end{abstract}

\noindent{\qquad\!\bf Keywords.} Riemann-Hilbert problems, Clifford algebra, biaxial symmetry, Dirac operator.
\smallbreak

\noindent{\qquad\!\bf AMS Subject Classifications.} 30G35; 15A66; 30E25; 35Q15; 31B10.

%%%%%%%%%%%%%%%%%%%%%%%%%%%%%%%%%%%%%%%%%%%%%%%%%%%%%%%%%%%%%%%%%%%%%%%%%%%%%%%%%%%%%%%%%%%%%%%%%%%%%%%%%%%%%%%%
%%%%%%%%%%%%%%%%%%%%%%%%%%%%%%%%%%%%%%%%%%%%%%%%%%%%%%%%%%%%%%%%%%%%%%%%%%%%%%%%%%%%%%%%%%%%%%%%%%%%%%%%%%%%%%%%%
\section{Introduction}\label{sec1}

\noindent

\noindent The Riemann-Hilbert boundary value problem (RHBVP), first formulated by Riemann in 1851 for analytic functions, originated as a geometric problem with profound implications in complex analysis \cite{R}. Hilbert’s seminal 1905 work reframed the problem using singular integral operators and the Hilbert transform \cite{H}, catalyzing a century of advancements by scholars including Muskhelishvili, Gakhov, and Lu  \cite{Ga,Lu,Mu}. These efforts expanded the theory to address diverse boundary value problems like Dirichlet, Neumann, Schwarz, and Robin problems on domains such as disks, half-planes, and annuli \cite{Ga,En,Mu1}. Further extensions encompassed null solutions of complex differential equations, including poly-analytic, meta-analytic, and poly-harmonic functions, bridging complex analysis with partial differential equations \cite{Fo,Be,Be1}.\\

\noindent The advent of Clifford analysis, boundary value problems including the Riemann-Hilbert problems, in higher-dimensional spaces are systematically addressed through its unified algebraic and geometric structure, transcending classical complex analysis. At the core of this framework lie monogenic functions, the null solutions to the Dirac or the generalized Cauchy-Riemann operators, which serve as the foundational tools for formulating and resolving multi-variable boundary value problems. In four dimensions, this theory aligns with Quaternionic analysis, providing a natural extension to hyper-complex systems. Early advancements established rigorous solutions to Riemann-Hilbert boundary value problems (RHBVPs) with constant coefficients for monogenic functions, leveraging integral representations and series expansions. While an exhaustive review of prior studies lies beyond the scope of this work, we highlight key advances in Clifford-analytic Riemann-Hilbert problems \cite{GD,ABP,Bu,GZ}. These results not only solidified the theoretical underpinnings of Clifford analysis but also bridged its applications to geometric configurations and operator-driven challenges in multi-variable domains. Building on this foundation, the progression from quaternionic to Clifford-analytic frameworks has been marked by pivotal advances in constructing monogenic functions. Fueter pioneered a method to generate monogenic functions from analytic counterparts in quaternionic spaces \cite{Fu}, a strategy later systematically extended by Sce to odd-dimensional Clifford algebras \cite{Sce} and by Qian to even-dimensional settings \cite{Qian}. These extensions catalyzed the study of variable coefficient boundary value problems, where non-constant operator coefficients introduced new analytical complexities \cite{He,He1,Ku,Ku1}. Subsequent studies further expanded these techniques to address axially symmetric configurations, resolving variable coefficient RHBVPs for monogenic \cite{Huang} and poly-monogenic functions \cite{Huang1}. Recent breakthroughs have generalized Fueter’s theorem to bi-axially symmetric domains, uncovering profound connections between generalized analytic functions and bi-axially symmetric monogenic solutions \cite{Pe,Ee,Ee2,Pe3,Pe2}. Specifically, two distinct classes of such functions were shown to satisfy the Vekua-type systems, analogous to generalized analytic function theories in higher dimensions. This interplay between symmetry, spectral perturbations, and hyper-complex operator theory now underpins modern approaches to multidimensional BVPs.\\

\noindent Despite transformative innovations in symmetry-driven Clifford analysis, the Riemann-Hilbert problems for the null solutions of perturbed iterated Dirac operators with vector wave numbers in $\mathbb{R}^{n}$ remain unexplored although such vector parameters are central to modeling anisotropic wave propagation in photonic crystals \cite{J1}, topological quantum states \cite{H2}, and elastic composites \cite{K3,G4}, where directional spectral coupling governs wave dynamics.
A gap motivated by their critical role in modeling multidimensional wave propagation, quantum fields, and elasticity with anisotropic material properties. Traditional methods for constant-coefficient or axially symmetric cases \cite{Huang1, Pe3} prove inadequate for these scenarios, where variable coefficients and vector-valued spectral parameters introduce nonlinear coupling and geometric complexity. In this work, we address this challenge by unifying bi-axial symmetry principles with a novel Almansi-type decomposition theorem, tailored to decouple poly-monogenic functions into hierarchically structured components. This decomposition enables explicit integral representations for solutions to variable-coefficient RHBVPs, even under vector wave number perturbations. Crucially, our method reveals that the Clifford-algebraic structure inherently parametrizes spectral distortions caused by wave number variations, transforming nonlinear interactions into algebraically tractable forms. As a cornerstone application, we resolve the Schwarz problem for these systems, obtaining closed-form solutions that generalize classical results to higher-dimensional, perturbed regimes. \\

\noindent This work establishes three pivotal advancements: (1) a theoretical framework for solving Riemann-Hilbert problems with Clifford-algebra-valued variable coefficients and vector wave number perturbations in bi-axially symmetric domains, generalizing higher-dimensional solvability (See Section $4$); (2) a novel Almansi-type decomposition methodology that enables direct computation of poly-monogenic solutions via recursive integral operators, overcoming spectral limitations of classical series expansions (See Section 3); and (3) a unified approach integrating classical complex-analytic tools into Clifford-theoretic calculus to solve heterogeneous boundary value problems across dimensions (See Section 4). By generalizing Fueter’s theorem and axial symmetry results while pioneering techniques for multi-variable spectral perturbations, this study bridges Clifford analysis with mathematical physics and spectral geometry, expanding its applicability to anisotropic and geometrically constrained systems.\\

\noindent The structure of this paper is designed to systematically advance from foundational concepts to novel applications. Section 2 presents essential preliminaries in Clifford analysis and rigorously formalizes bi-axially symmetric monogenic functions, elucidating their intrinsic relationship with Vekua-type governing equations. Section 3 develops a suite of instrumental lemmas, which underpin the subsequent analytical framework. In Section 4, we resolve Riemann-Hilbert boundary value problems for bi-axially symmetric poly-monogenic functions, the null solutions to iterated Dirac operators by constructing explicit solutions through a tailored Almansi-type decomposition theorem. Further, we generalize these results to perturbed iterated Dirac systems with vector wave number parameters, demonstrating the adaptability of our approach to spectral variations. As a culminating application, we derive closed-form solutions to the Schwarz problem, unifying classical techniques with Clifford-algebraic operator calculus. This progression highlights how symmetry-driven decomposition strategies transcend dimensional constraints, offering a blueprint for higher-order boundary value challenges in multi-variable domains.

\medskip

\section{Preliminaries}
In this section, we briefly review some notions and definitions of Clifford algebra. We refer to \cite{De,Lo} for more details.\\

\noindent Let $\{e_1,\ldots,e_n\},n\in\mathbb{N}$ be the $n$ orthonormal basis of the Euclidean space $\mathbb{R}^{n}$ and the multiplication satisfies the following relations:
$e^{2}_{j}=-1$, if $j=1,2,\ldots,n$ and $e_{i}e_{j}+e_{j}e_{i}=0$, if $1\leq i\neq j\leq n$.
The conjugation $\bar{e}_j$ of $e_j$ is defined by $\bar{e}_j=-e_j, j=1,2,\ldots,n$. Thus, any Clifford number $a$ in real Clifford algebra $\mathbb{R}_{0,n}$ can be written as
\begin{eqnarray*}
a=\sum_{A}a_{A}e_{A}, \ \ a_{A}\in\mathbb{R},
\end{eqnarray*}
where $e_{A}=e_{j_1}e_{j_2} \ldots e_{j_k}$, $A=\{j_1,j_2,\ldots,j_k\}\subseteq\{1,2,\ldots,n\}$ with $j_1<\ldots<j_k$.
Especially, for $A=\varnothing$, we write $e_\varnothing$ as $e_0$ and put $e_\varnothing=1$.
The inner product of $a,b\in \mathbb{R}_{0,n}$ is defined by
\begin{eqnarray}\label{2.1}
\left(a,b\right) = \left[ab\right]_0,
\end{eqnarray}
where $\left[a\right]_0$ denotes the scalar part of $a$.\\

\noindent Furthermore, Euclidean space $\mathbb{R}^{n}$ can be embedded in Clifford algebra $\mathbb{R}_{0,n}$ through the element $(x_1,\ldots,x_n)\in \mathbb{R}^{n}$ and the vector $\underline{x}=\sum^n_{j=1} x_je_j$. The conjugation $\overline{\underline{x}}$ of the vector $\underline{x}$ is written as $\overline{\underline{x}}=\sum^n_{j=1} x_j\bar{e}_j=-\underline{x}$.
The Euclidean norm of arbitrary $\underline{x}\in\mathbb{R}^{n}$ is $|\underline{x}|=(\sum\nolimits^{n}_{j=1}x_{j}^{2})^{\frac{1}{2}}$, and the multiplication of $\underline{x}$ satisfies $\underline{x}^2=-|\underline{x}|^2=-\sum^{n}_{j=1}x_j^2$. If $\underline{x}\in\mathbb{R}^{n}\backslash\{0\}$, then the inverse $\underline{x}^{-1}$ exists and $\underline{x}^{-1}:=\overline{\underline{x}}\cdot |\underline{x}|^{-2}$, i.e. $\underline{x}\underline{x}^{-1}=\underline{x}^{-1}\underline{x}=1$.\\

\noindent We investigate the bi-axial decomposition $\mathbb{R}^n=\mathbb{R}^p\oplus\mathbb{R}^q$ with $p+q=n$ (see \cite{Pe2,Qi}) by expanding Fueter's theorem \cite{Fu} to the biaxial scenario. Hence, we write any $\underline x\in\mathbb{R}^{n}$ as
\begin{eqnarray*}
\underline{x}=\underline{y}+\underline{z}
\end{eqnarray*}
with
\begin{eqnarray*}
\underline{y}=\sum_{j=1}^{p}y_{j}e_{j},\underline{z}=\sum_{j=p+1}^{n}z_{j}e_{j}.\\
\end{eqnarray*}

\noindent We now define the bi-axially symmetric open sets, referring to the references \cite{GJ,AM}.\\

\begin{De}
Let $\Omega$ be a non-empty open subset in $\mathbb{R}^{n}$. For arbitrary $\underline{x}\in\Omega$, we define the set
\begin{eqnarray}\label{2.2}
[\underline{x}]=\{\underline{x}:\ \underline{x}=P\underline{y}+Q\underline{z},\ P\in \mathrm{SO}(p),\  Q\in \mathrm{SO}(q)\},
\end{eqnarray}
where $\mathrm{SO}(p)$ and $\mathrm{SO}(q)$ are the groups containing standard orthogonal transformations. We call that $\Omega$ is bi-axially symmetric if the set $[\underline{x}]$ is contained in $\Omega\subset{\mathbb{R}^{n}}$ for arbitrary $\underline{x}\in\Omega$, i.e. $\Omega$ is invariant under $\mathrm{SO}(p)\times \mathrm{SO}(q)$.
\end{De}

\begin{Rem}
The ball centered at the origin of $\mathbb{R}^{n}$ is a bi-axially symmetric domain.
\end{Rem}

\noindent The first order differential operator in $\mathbb{R}^n$ written as
\begin{eqnarray}
\mathcal{D}=\sum_{j=1}^{n}e_j\partial_{x_j}
\end{eqnarray}
is called the Dirac operator, and monogenic functions are functions in the kernel of $\mathcal{D}$\cite{De}. It factorizes the Laplacian $\Delta= - \mathcal{D}^2 =\sum^{n}_{j=1}\partial^{2}_{x_j}$ in $\mathbb{R}^n$.\\

\noindent Introduce the perturbed Dirac operator with a vector wave number in $\mathbb{R}^n$
\begin{eqnarray}
\mathcal{D} - \alpha =\sum_{j=1}^{n}e_j\left(\partial_{x_j} - \alpha_j\right),
\end{eqnarray}
where $\alpha=\sum\limits_{j=1}^{n}e_{j}\alpha_{j}\in\mathbb{R}^{n}, \alpha_{j} \in \mathbb{R}$. Then, it factorizes the following second-order elliptic differential operator in $\mathbb{R}^n$, that is,
\begin{eqnarray}
    \left(\mathcal{D} - \alpha\right)^2 = -\left(\Delta + \sum\limits^n_{j=1}\alpha_{j}^2 \right)- 2\sum\limits^{n}_{j=1} \alpha_{j} \partial_{x_j} + \sum\limits_{1\leq i < j \leq n}e_i e_j\left(\alpha_{i} \partial_{x_j} - \alpha_{j} \partial_{x_i}\right).
\end{eqnarray}

\noindent Let $\underline{x}=\underline{y}+\underline{z}\in\mathbb{R}^{n}, r=|\underline{y}|,\rho=|\underline{z}|$, $\underline{\omega}\in S^{p-1}\subset\mathbb{R}^p,  \underline{\nu}\in S^{q-1}\subset\mathbb{R}^q $.

\noindent Then, the Dirac operator in spherical coordinates on the bi-axial splitting can be expressed as
\begin{eqnarray*}
\mathcal{D}=\underline{\omega}(\partial_r+\frac{1}{r}\Gamma_{\underline{\omega}})+\underline{\nu}(\partial_{\rho}+\frac{1}{\rho}\Gamma_{\underline{\nu}})
\end{eqnarray*}
with $\Gamma_{\underline{\omega}}=-\underline{\omega}\wedge\mathcal{D}_{\underline{\omega}}=\overline{\underline{\omega}}\mathcal{D}_{\underline{\omega}},
\Gamma_{\underline{\nu}}=-\underline{\nu}\wedge\mathcal{D}_{\underline{\nu}}=\overline{\underline{\nu}}\mathcal{D}_{\underline{\nu}}$
and $\mathcal{D}_{\underline{\omega}}=\sum^{p}_{j=1}e_{j}\partial_{y_j},\mathcal{D}_{ \underline{\nu}}=\sum^{n}_{j=p+1}e_{j}\partial_{z_j}$.\\

\noindent In this paper, let $\Omega\subset\mathbb{R}^n$ be a bi-axially symmetric domain with a smooth boundary $\partial\Omega$. The function $\varphi:\Omega\rightarrow\mathbb{R}_{0,n}$ is denoted by $\varphi(\underline{x})=\sum\varphi_A(\underline x)e_A$, where $\varphi_A(\underline x)$ is a $\mathbb{R}$-valued function. The function $\varphi$ is continuous, continuously differentiable, H\"{o}lder continuous means that each component of $\varphi$ has the same properties. The corresponding spaces are denoted by $\mathcal{C}(\Omega, \mathbb{R}_{0,n}), \mathcal{C}^1(\Omega, \mathbb{R}_{0,n}), \mathcal{H}^{\mu}(\Omega, \mathbb{R}_{0,n})$ $(0 <\mu \leq 1)$, respectively.\\

\begin{De}
If $\varphi\in\mathcal{C}^1(\Omega, \mathbb{R}_{0,n})$, satisfying $\mathcal{D}\varphi(\underline x)=0$, it is called (left)-monogenic function.
\end{De}

\noindent If $\varphi(\underline{x})$ is of the form
\begin{eqnarray}\label{2.3}
\varphi(\underline{x})=A(r,\rho)+\underline{\omega}\underline{\nu}B(r,\rho)+\underline{\omega}C(r,\rho)+\underline{\nu}D(r,\rho),
\end{eqnarray}
where $A,B,C,D$ are $\mathbb{R}$-valued continuously differentiable functions in $\mathbb{R}^2$, we call it of bi-axial type.\\

\begin{De}
A function $\varphi(\underline{x})\in\mathcal{C}^1(\Omega, \mathbb{R}_{0,n})$ is said to be bi-axially monogenic, if it is monogenic and of bi-axial type.
\end{De}

\begin{Rem}\label{Rem2.2}
The expression $\mathcal{D}(A+\underline{\omega}\underline{\nu}B+\underline{\omega}C+\underline{\nu}D)$ is of the form $a\underline{\omega}+b\underline{\nu}+c+\underline{\omega}\underline{\nu}d$, where $a,b,c,d$ are $\mathbb{R}$-valued functions satisfying $a\underline{\omega}+b\underline{\nu}=\mathcal{D}(A+\underline{\omega}\underline{\nu}B),
c+\underline{\omega}\underline{\nu}d=\mathcal{D}(\underline{\omega}C+\underline{\nu}D)$.\\

\noindent Furthermore, with respect to the inner product \eqref {2.1}, $a\underline{\omega}+b\underline{\nu}$ is orthogonal to $c+\underline{\omega}\underline{\nu}d$. \\

\noindent So, the equation $\mathcal{D}(A+\underline{\omega}\underline{\nu}B+\underline{\omega}C+\underline{\nu}D)=0$ can be split into two independent equations
\begin{eqnarray}
\left\{\begin{array}{ll}
\mathcal{D}(A+\underline{\omega}\underline{\nu}B)=0,\\
\mathcal{D}(\underline{\omega}C+\underline{\nu}D)=0.
\end{array}\right.
\end{eqnarray}
Hence, we obtain two special types of bi-axially monogenic functions
\begin{eqnarray}\label{2.4}
\varphi_1(\underline{x})=A(r,\rho)+\underline{\omega}\underline{\nu}B(r,\rho),
\end{eqnarray}
and
\begin{eqnarray}\label{2.5}
\varphi_2(\underline{x})=\underline{\omega}C(r,\rho)+\underline{\nu}D(r,\rho).
\end{eqnarray}
We call $\varphi_1$ and $\varphi_2$ as the first type and the second type of bi-axially symmetric functions, respectively.
\end{Rem}

\noindent Next, we consider the first type of bi-axially monogenic functions. \\

\noindent As $\Gamma_{\underline{\omega}}A=0,\Gamma_{\underline{\nu}}A=0,
\Gamma_{\underline{\omega}}(\underline{\omega}\underline{\nu}B)=(p-1)\underline{\omega}\underline{\nu}B,
\Gamma_{\underline{\nu}}(\underline{\omega}\underline{\nu}B)=(q-1)\underline{\omega}\underline{\nu}B$, we have
\begin{eqnarray*}
0=\mathcal{D}\varphi_1&=&\underline{\omega}\partial_{r}A+\underline{\omega}\partial_{r}(\underline{\omega} \underline{\nu}B)+\frac{\underline{\omega}}{r}\Gamma_{\underline{\omega}}A+\frac{\underline{\omega}}{r}\Gamma_{\underline{\omega}}(\underline{\omega} \underline{\nu}B)  \\
&\ \ \ +&\underline{\nu}\partial_{\rho}A+\underline{\nu}\partial_{\rho}(\underline{\omega} \underline{\nu}B)+\frac{\underline{\nu}}{\rho}\Gamma_{\underline{\nu}}A+\frac{\underline{\nu}}{\rho}\Gamma_{\underline{\nu}}(\underline{\omega} \underline{\nu}B)  \\
&=&\underline{\omega}\partial_{r}A-\underline{\nu}\partial_{r}B-\frac{p-1}{r}\underline{\nu}B+\underline{\nu}\partial_{\rho}A+
\underline{\omega}\partial_{\rho}B+\frac{q-1}{\rho}\underline{\omega}B  \\
&=&\underline{\omega}(\partial_rA+\partial_{\rho}B+\frac{q-1}{\rho}B)+\underline{\nu}(\partial_{\rho}A-\partial_{r}B-\frac{p-1}{r}B).  \\
\end{eqnarray*}
It is clear that the functions $A,B$ satisfy the Vekua-type system
\begin{eqnarray}\label{2.6}
\left\{\begin{array}{ll} \partial_{r}A+\partial_{\rho}B=-\frac{q-1}{\rho}B,\\
\partial_{\rho}A-\partial_{r}B=\frac{p-1}{r}B.
\end{array}\right.
\end{eqnarray}

\noindent With regard to the second type of bi-axially monogenic functions, because of
$$\Gamma_{\underline{\nu}}(\underline{\omega}C)=0, \Gamma_{\underline{\omega}}(\underline{\nu}D)=0,
\Gamma_{\underline{\omega}}(\underline{\omega}C)=(p-1)\underline{\omega}C,\Gamma_{\underline{\nu}}(\underline{\nu}D)=(q-1)\underline{\nu}D,$$
we have
\begin{eqnarray*}
0=\mathcal{D}\varphi_2&=&\underline{\omega}\partial_{r}(\underline{\omega}C)+\underline{\omega}\partial_{r}(\underline{\nu}D)
+\frac{\underline{\omega}}{r}\Gamma_{\underline{\omega}}(\underline{\omega}C)+\frac{\underline{\omega}}{r}\Gamma_{\underline{\omega}}(\underline{\nu}D)  \\
&\ \ \ +&\underline{\nu}\partial_{\rho}(\underline{\omega}C)+\underline{\nu}\partial_{\rho}(\underline{\nu}D)
+\frac{\underline{\nu}}{\rho}\Gamma_{\underline{\nu}}(\underline{\omega}C)+\frac{\underline{\nu}}{\rho}\Gamma_{\underline{\nu}}(\underline{\nu}D)  \\
&=&-\partial_{r}C+\underline{\omega}\underline{\nu}\partial_{r}D-\frac{p-1}{r}C-\underline{\omega}\underline{\nu}\partial_{\rho}C-
\partial_{\rho}D-\frac{q-1}{\rho}D  \\
&=&-(\partial_rC+\partial_{\rho}D+\frac{p-1}{r}C+\frac{q-1}{\rho}D)-\underline{\omega}\underline{\nu}(\partial_{\rho}C-\partial_{r}D).
\end{eqnarray*}

\noindent Then, the functions $C,D$ satisfy another Vekua-type system
\begin{eqnarray}\label{2.7}
\left\{\begin{array}{ll} \partial_{r}C+\partial_{\rho}D=-\frac{p-1}{r}C-\frac{q-1}{\rho}D,\\
\partial_{\rho}C-\partial_{r}D=0.
\end{array}\right.
\end{eqnarray}

\begin{De}
A function $\varphi(\underline{x})\in \mathcal{C}^k(\Omega, \mathbb{R}_{0,n}),k\in\mathbb{N}$ is said to be a bi-axially (left) $k$-monogenic function or poly-monogenic function, if it is of biaxial type and satisfies $\mathcal{D}^k\varphi(\underline{x})=0$, where $\mathcal{D}^{k} \varphi \triangleq \mathcal{D}^{k-1}(\mathcal{D} \varphi)$.
\end{De}

\noindent In the following, any bi-axially symmetric function is defined on $\overline{\Omega}\subset\mathbb{R}^n$ where $\Omega$ is bi-axially symmetric
without explanation. $S\subset\mathbb{C}_+$ is the projection of the biaxially symmetric domain $\Omega\subset\mathbb{R}^n$ into the $(r,\rho)$-plane, where $\mathbb{C}_+$ is the quarter of the $(r,\rho)$-plane.
\medskip

\section{Some lemmas}
In this section, We present several lemmas concerning the Euler operator and introduce a Almansi-type decomposition for null solutions of the iterated perturbed Dirac operator with a vector wave number, which will be utilized subsequently.

\begin{De}\label{de3.1}
Let $\Omega\subset\mathbb{R}^n$ be a star-like subdomain with center $0$ if with $\underline{x}\in\Omega$ also $t\underline{x}\in\Omega$ holds for $\forall 0\leq t \leq 1$. The Euler operator defined on the space $\mathcal{C}^{1}\left(\Omega, \mathbb{R}^{n}\right)$ is given by
\begin{eqnarray}\label{3.1}
E_{s}=s I+\sum_{j=1}^{n}x_{j} \partial_{x_{j}},\ s>0,
\end{eqnarray}
where $I$ denotes the identity operator.
\end{De}

\begin{De}\label{de3.2}
Let $\Omega\subset\mathbb{R}^n$ be a star-like sub-domian with the center at $0$. The operator $I_{s}: \mathcal{C}\left(\Omega, \mathbb{R}^{n}\right) \rightarrow \mathcal{C}\left(\Omega, \mathbb{R}^{n}\right)$
is defined by
\begin{eqnarray}\label{3.2}
I_{s} \varphi=\int_{0}^{1} \varphi(t \underline{x}) t^{s-1} d t,\ s>0, \underline{x} \in \Omega.
\end{eqnarray}

\noindent For arbitrary $k\in\mathbb{N}$, we denote
\begin{eqnarray}\label{3.3}
\mathcal{Q}_{k}=\frac{1}{a_{k}} I_{\frac{n}{2}} I_{\frac{n}{2}+1} \ldots I_{\frac{n}{2}+\left[\frac{k-1}{2}\right]},\ a_{k}=(-2)^{k}\left[\frac{k}{2}\right]!,
\end{eqnarray}
where
$$
[s]= \begin{cases}k, & \text { if } k \in \mathbb{N}, \\
k+1, & \text { if } s=k+t, k \in \mathbb{N}, 0<t<1 .\end{cases}
$$
\end{De}

\begin{lem} \cite{Ku2} \label{lem3.1}
Suppose $\Omega$ is a star-like sub-domain with the center at $0$, we have
\begin{eqnarray}\label{3.4}
\varphi=E_{s} I_{s}\varphi=I_{s} E_{s}\varphi,
\end{eqnarray}
and
\begin{eqnarray}\label{3.5}
\mathcal{D} E_{s} \varphi=E_{s+1} \mathcal{D} \varphi.
\end{eqnarray}
\end{lem}

\begin{Rem}\label{Rem3.1}
If $\varphi\in(\Omega,\mathbb{R}^n)$ is monogenic, we have that $E_{s} \varphi$ and $I_{s} \varphi$ are both monogenic by \eqref{3.5}.

\noindent Indeed, when $\mathcal{D}^2 \varphi=0$, we have
$$
\mathcal{D}^2(E_s \varphi)=\mathcal(E_{s+1}\mathcal{D}\varphi)=E_{s+2}\mathcal{D}^2 \varphi.
$$
Then, we can obtain
$$
\mathcal{D}^k(E_s \varphi)=E_{s+k}\mathcal{D}^k \varphi,\ k\in\mathbb{N},
$$
in an inductive way. \\

\noindent It is also similar to the operator $I_s$. \\

\noindent Therefore, if $\varphi\in(\Omega,\mathbb{R}^n)$ is $k$-monogenic, we have that $E_{s} \varphi$ and $I_{s} \varphi$ are both $k$-monogenic.\\

\noindent Moreover, suppose $\Omega$ is of bi-axial type and is a star-like subdomain with the center at $0$. By the direct calculation, we have
\begin{eqnarray*}
E_s \varphi(\underline{x}) &=& s\varphi(\underline{x})+ \sum_{j=1}^{p}x_{j}  \frac{\partial \varphi}{\partial x_{j}}(\underline{x})+\sum_{j=p+1}^{n}x_{j} \frac{\partial \varphi}{\partial x_{j}}(\underline{x}) \\
&=& s\varphi(\underline{x})+ (r\partial r+\rho \partial \rho) \varphi(\underline{x}),\underline{x}\in\Omega.
\end{eqnarray*}

\noindent For $I_s$, the integral is made for $t$ and does not change the bi-axial symmetry of $\varphi$. Hence, if $\varphi\in(\Omega,\mathbb{R}^n)$ is of bi-axial type, we have that $E_{s} \varphi$ and $I_{s} \varphi$ are both of bi-axial type.
\end{Rem}

\begin{lem} \label{lem3.2}
Let $\Omega\subset\mathbb{R}^n$ be a star-like sub-domain with the  center at $0$. If $\varphi \in\mathcal{C}^{k}\left(\Omega, \mathbb{R}_{0,n}\right)$ is a solution to $(\mathcal{D}-\alpha)^{k} \varphi=0, \alpha=\sum_{j=1}^{n}e_{j}\alpha_{j}\in\mathbb{R}^{n}$, then there exist uniquely defined functions $\varphi_{j}$, such that
\begin{eqnarray}\label{6}
\varphi(\underline{x})=\mathrm{e}^{\beta} \sum_{j=0}^{k-1} \underline{x}^{j} \varphi_{j}(\underline{x}), \underline{x} \in \Omega,
\end{eqnarray}
where $\beta=\sum_{j=1}^{n}\alpha_{j}x_{j}$, each $\varphi_j$ is monogenic and given by
\begin{eqnarray}\label{8}
\left\{\begin{array}{l}
\varphi_{0}=\left(I-\underline{x} \mathcal{Q}_{1} \mathcal{D}\right)\left(I-\underline{x}^{2} \mathcal{Q}_{2} \mathcal{D}^{2}\right) \ldots\left(I-\underline{x}^{k-1} \mathcal{Q}_{k-1} \mathcal{D}^{k-1}\right) (\mathrm{e}^{-\beta}\varphi), \\
\varphi_{1}=\mathcal{Q}_{1} \mathcal{D}\left(I-\underline{x}^{2} \mathcal{Q}_{2} \mathcal{D}^{2}\right) \ldots\left(I-\underline{x}^{k-1} \mathcal{Q}_{k-1} \mathcal{D}^{k-1}\right) (\mathrm{e}^{-\beta}\varphi), \\
\quad \vdots \\
\varphi_{k-2}=\mathcal{Q}_{k-2} \mathcal{D}^{k-2}\left(I-\underline{x}^{k-1} \mathcal{Q}_{k-1} \mathcal{D}^{k-1}\right) (\mathrm{e}^{-\beta}\varphi), \\
\varphi_{k-1}=\mathcal{Q}_{k-1} \mathcal{D}^{k-1} (\mathrm{e}^{-\beta}\varphi).
\end{array}\right.
\end{eqnarray}
\end{lem}

\noindent\p\  Since $\alpha\in\mathcal{C}^{k}\big{(}\Omega,\mathbb{R}_{0,n}\big{)},\ \mathcal{D}\beta=\alpha$, we have
\begin{eqnarray*}
\big{(}\mathcal{D}-\alpha \big{)}\varphi(\underline{x})= \big{(}\mathcal{D}-\alpha \big{)} \mathrm{e}^{\beta} \big{(} \mathrm{e}^{-\beta}\varphi(\underline{x})\big{)}=\mathrm{e}^{\beta} \mathcal{D}\big{(} \mathrm{e}^{-\beta}\varphi(\underline{x})\big{)}.
\end{eqnarray*}
By directly calculating, we get
\begin{eqnarray*}
\big{(}\mathcal{D}-\alpha\big{)}^{k}\varphi(\underline{x})=\mathrm{e}^{\beta} \mathcal{D}^{k}\big{(} \mathrm{e}^{-\beta}\varphi(\underline{x})\big{)},
\end{eqnarray*}
which means that the function $\mathrm{e}^{-\beta}\varphi(\underline{x})$ is $k$-monogenic.
From Theorem 2.1 in \cite{Ma}, we can obtain
$$
\mathrm{e}^{- \beta} \varphi(\underline{x})=\sum_{j=0}^{k-1} \underline{x}^{j} \varphi_{j}(\underline{x}), \underline{x} \in \Omega,
$$
where each $\varphi_j$ is uniquely given by expression
\begin{eqnarray*}
\left\{\begin{array}{l}
\varphi_{0}=\left(I-\underline{x} \mathcal{Q}_{1} \mathcal{D}\right)\left(I-\underline{x}^{2} \mathcal{Q}_{2} \mathcal{D}^{2}\right) \ldots\left(I-\underline{x}^{k-1} \mathcal{Q}_{k-1} \mathcal{D}^{k-1}\right) (\mathrm{e}^{-\beta}\varphi), \\
\varphi_{1}=\mathcal{Q}_{1} \mathcal{D}\left(I-\underline{x}^{2} \mathcal{Q}_{2} \mathcal{D}^{2}\right) \ldots\left(I-\underline{x}^{k-1} \mathcal{Q}_{k-1} \mathcal{D}^{k-1}\right) (\mathrm{e}^{-\beta}\varphi), \\
\quad \vdots \\
\varphi_{k-2}=\mathcal{Q}_{k-2} \mathcal{D}^{k-2}\left(I-\underline{x}^{k-1} \mathcal{Q}_{k-1} \mathcal{D}^{k-1}\right) (\mathrm{e}^{-\beta}\varphi), \\
\varphi_{k-1}=\mathcal{Q}_{k-1} \mathcal{D}^{k-1} (\mathrm{e}^{-\beta}\varphi),
\end{array}\right.
\end{eqnarray*}
and each $\varphi_j$ is monogenic. Hence, we get the result
\begin{eqnarray*}
\varphi(\underline{x})=\mathrm{e}^{\beta} \sum_{j=0}^{k-1} \underline{x}^{j} \varphi_{j}(\underline{x}), \underline{x} \in \Omega.
\end{eqnarray*}
$\hfill\square$

\begin{Rem}
    Let $Ker (\mathcal{D}-\alpha)^{k} = \big{\{} \phi: (\mathcal{D}-\alpha)^{k}\phi =0 \big{\}}, k\in \mathbb{N}, \alpha=\sum\limits_{j=1}^{n}e_{j}\alpha_{j}\in\mathbb{R}^{n}$. Then, one has
    \begin{eqnarray}
        Ker (\mathcal{D}-\alpha)^{k} = \mathrm{e}^{-\beta} \bigoplus^{k-1}_{j=0} \underline{x}^{j} Ker (\mathcal{D}-\alpha),
    \end{eqnarray}
  where $\beta=\sum\limits_{j=1}^{n}\alpha_{j}x_{j}$.
  \noindent
\end{Rem}

\begin{Rem}
    Lemma $\ref{lem3.2}$ gives the direct sum decomposition of the for null-solutions to iterated perturbed Dirac operator with a vector wave number, which is called new Almansi-type decomposition (in short Almansi-type decomposition). Similar Almansi-type decomposition for null-solutions to iterated perturbed Dirac operator with a scalar wave number Refers to Refs. e.g. \cite{Ku2,Ma}.
\end{Rem}

\begin{Cor}\label{Cor3.1}
Let $\Omega\subset\mathbb{R}^n$ be a star-like sub-domain with the center at $0$. If $\varphi \in\mathcal{C}^{k}\left(\Omega, \mathbb{R}_{0,n}\right)$ is a solution to $\mathcal{D}^{k} \varphi=0$, then there exist uniquely defined functions $\varphi_{j}$ such that
\begin{eqnarray}\label{3.6}
\varphi=\sum_{j=0}^{k-1} \underline{x}^{j} \varphi_{j}, \underline{x} \in \Omega,
\end{eqnarray}
where each $\varphi_{j}(j=0,1,2, \ldots, k-1)$ is monogenic in $\Omega$ and given by
\begin{eqnarray}\label{3.7}
\left\{\begin{array}{l}
\varphi_{0}=\left(I-\underline{x} \mathcal{Q}_{1} \mathcal{D}\right)\left(I-\underline{x}^{2} \mathcal{Q}_{2} \mathcal{D}^{2}\right) \ldots\left(I-\underline{x}^{k-1} \mathcal{Q}_{k-1} \mathcal{D}^{k-1}\right) \varphi, \\
\varphi_{1}=\mathcal{Q}_{1} \mathcal{D}\left(I-\underline{x}^{2} \mathcal{Q}_{2} \mathcal{D}^{2}\right) \ldots\left(I-\underline{x}^{k-1} \mathcal{Q}_{k-1} \mathcal{D}^{k-1}\right) \varphi, \\
\quad \vdots \\
\varphi_{k-2}=\mathcal{Q}_{k-2} \mathcal{D}^{k-2}\left(I-\underline{x}^{k-1} \mathcal{Q}_{k-1} \mathcal{D}^{k-1}\right) \varphi, \\
\varphi_{k-1}=\mathcal{Q}_{k-1} \mathcal{D}^{k-1} \varphi.
\end{array}\right.
\end{eqnarray}
\end{Cor}

\noindent By Remark \ref{2.2}, Remark \ref{Rem3.1} and Lemma \ref{lem3.2}, we have the following corollary.

\begin{Cor}
Suppose $\Omega$ is of bi-axial type and is a star-like sub-domain with the center at $0$. If $\varphi \in\mathcal{C}^{k}\left(\Omega, \mathbb{R}_{0,n}\right)$ satisfying $(\mathcal{D}-\alpha)^{k} \varphi=0$ is of biaxial type, then each $\varphi_j(j=0,1,2, \ldots, k-1)$ given by equality \eqref{6} is of biaxial type.\\

\noindent If $\varphi \in\mathcal{C}^{k}\left(\Omega, \mathbb{R}_{0,n}\right)$ satisfying $\mathcal{D}^{k} \varphi=0$ is of biaxial type, then each $\varphi_j(j=0,1,2, \ldots, k-1)$ given by equality \eqref{3.7} is of biaxial type.
\end{Cor}

\noindent\p\ We know that the sum and product of two bi-axial functions are still of biaxial type by direct calculation. For $k-1$, we can obtain $\mathcal{D}(\mathrm{e}^{-\beta}\varphi)$ is of biaxial type by Remark \ref{Rem2.2}. Repeating $k-2$ times, we derive $\mathcal{D}^{k-1}(\mathrm{e}^{-\beta}\varphi)$ is of biaxial type. On the other hand, as biaxial symmetry holds after acting $I_s$, we can obtain that biaxial symmetry also holds after acting $\mathcal{Q}_{k-1}$. Therefore, function $\varphi_{k-1}$ is a biaxial function.\\

\noindent For $k-2$, $\left(I-\underline{x}^{k-1} \mathcal{Q}_{k-1} \mathcal{D}^{k-1}\right)\varphi$ is a bi-axial function. Indeed, $\underline{x}$ can be written as $\underline{x}=\underline{\omega}r+\underline{\nu}\rho$, which can be seen as a biaxial function. Then, acting $\mathcal{D}^{k-1}$ and $\mathcal{Q}^{k-1}$, we obtain that function $\varphi_{k-2}$ is of biaxial type.\\

\noindent Hence, every $\varphi_j(j=0,1,\ldots,k-1)$ in equality \eqref{6} is of biaxial type in inductive way. Of course, every $\varphi_j(j=0,1,\ldots,k-1)$ in equality \eqref{3.7} is also of bi-axial type in the same way.
$\hfill\square$

\begin{lem}\cite{Ku2} \label{lem3.3}
Suppose that $\Omega$ is a sub-domain of $\mathbb{R}^{n}$ and $j \in \mathbb{N}$ is arbitrary. If $\varphi \in \mathcal{C}^{1}$  $\left(\Omega, \mathbb{R}_{0,n}\right)$ is monogenic, then
$$
\mathcal{D}\left(\underline{x}^{j} \varphi\right)= \begin{cases}-2 m \underline{x}^{2 m-1} \varphi, & \ if \ j=2 m, \\
-2 \underline{x}^{2(m-1)} E_{\frac{n+1}{2}+\left[\frac{j}{2}\right]-1} \varphi, & \ if \ j=2 m-1,\end{cases}
$$
with $m \in \mathbb{N}$. Moreover, for $l \in \mathbb{N}$ and $2 \leq l \leq j$, we get
\begin{eqnarray}
\mathcal{D}^{l} \underline{x}^{j} \varphi=\mathfrak{C}_{l, j} \underline{x}^{j-l} E_{\frac{n+1}{2}+\left[\frac{j-l}{2}\right]} \ldots E_{\frac{n+1}{2}+\left[\frac{j}{2}\right]-1} \varphi ,
\end{eqnarray}
where
$$
\mathfrak{C}_{l, j}= \begin{cases}2^{l} m(m-1) \ldots(m-p+1), & \text { if } j=2 m, l=2 p, \\
-2^{l} m(m-1) \ldots(m-p), & \ if \ j=2 m, l=2 p+1, \\ 2^{l}(m-1) \ldots(m-p), & \ if \ j=2 m-1, l=2 p, \\
-2^{l}(m-1) \ldots(m-p+1), & \ if \ j=2 m-1, l=2 p-1 ,\end{cases}
$$
with $m,p \in \mathbb{N}.$
\end{lem}

\begin{Rem}\label{Rem3.3}
For $l=j$, we obtain
$$
\mathcal{D}^{j} \underline{x}^{j} \varphi=\mathfrak{C}_{j,j} E_{\frac{n+1}{2}} \ldots E_{\frac{n+1}{2}+m-1} \varphi
$$
with
$$\mathfrak{C}_{j, j}=\begin{cases}2^{j} m!, &\ if \ j=2 m, \\
-2^{j}(m-1)!, &\ if \ j=2 m-1 . \end{cases}
$$
\end{Rem}

\begin{lem}\cite{Zuo} \label{lem3.4}
The RHBVP with Clifford-valued variable coefficients for monogenic functions of bi-axial type is: to find a function
$\varphi\in\mathcal{C}^{1}\big(\Omega,\mathbb{R}_{0,n}\big)$, satisfying
\begin{equation}\label{3.9}
\left\{\begin{array}{ll}\mathcal{D}\varphi(\underline{x})=0,\ \ \underline{x}\in\Omega,\\
 \mathrm{Re}\{\lambda_1(\underline{t})\varphi(\underline{t})\}= g_1(\underline{t}),\ \ \underline{t}\in\partial\Omega,\\
 \mathrm{Re}\{\lambda_2(\underline{t})\varphi(\underline{t})\}= g_2(\underline{t}),\ \ \underline{t}\in\partial\Omega,
\end{array}\right.
\end{equation}
where $\lambda_1(\underline{t})=A_b(r,\rho)+\underline{\omega}\underline{\nu}B_b(r,\rho)$ is of the first bi-axial type, $\lambda_2(\underline{t})
=\underline{\omega}C_b(r,\rho)+\underline{\nu}D_b(r,\rho)$ is of the second biaxial type, and $\lambda_{i}\in\mathcal{H}^{\mu}(\partial\Omega,\mathbb{R}_{0,n}),g_{i}\in\mathcal{H}^{\mu}(\partial \Omega,\mathbb{R})(i=1,2)$.
$\mathrm{Re}\{\cdot\}$ can be seen as $[\cdot]_0$ in \eqref{2.1} and also in the following without explanation.\\

\noindent The RHBVP is solvable and its solution is given by
\begin{eqnarray}\label{3.10}
\varphi(\underline{x})=\mathrm{Im}f_{1}(r,\rho)+\underline{\omega\nu}\mathrm{Re}f_{1}(r,\rho)
+\underline{\omega}\mathrm{Im}f_{2}(r,\rho)+\underline{\nu}\mathrm{Re}f_{2}(r,\rho),\underline{x}\in\Omega,
\end{eqnarray}
where $\mathrm{Re}\{\cdot\},\mathrm{Im}\{\cdot\}$ denote the real part and imaginary part of complex valued functions $f_{i}(i=1,2)$ in $S\subset\mathbb{C}_{+}$ respectively. And $f_{i}$ are given by
\begin{eqnarray}\label{3.11}
f_{i}(z)=\Phi_{i}(z)\mathrm{e}^{\tau_{i}(z)},z=r+\mathrm{i}\rho,
\end{eqnarray}
with
\begin{eqnarray}\label{3.12}
\tau_{1}(z)=\frac1\pi\iint\limits_S\left((\frac{p-1}{4\xi}+\frac{q-1}{4\eta}\mathrm{i})+(\frac{p-1}{4\xi}+\frac{q-1}{4\eta}\mathrm{i})
\frac{\bar{f}_{1}(\zeta)}{f_{1}(\zeta)}\right)  \frac{\mathrm{d}\xi\mathrm{d}\eta}{\zeta-z},
\end{eqnarray}
and
\begin{eqnarray}\label{3.13}
\tau_{2}(z)=\frac1\pi\iint\limits_S\left((\frac{p-1}{4\xi}+\frac{q-1}{4\eta}\mathrm{i})-(\frac{p-1}{4\xi}-\frac{q-1}{4\eta}\mathrm{i})
\frac{\bar{f}_{2}(\zeta)}{f_{2}(\zeta)}\right)  \frac{\mathrm{d}\xi\mathrm{d}\eta}{\zeta-z},
\end{eqnarray}
where $\zeta=\xi+\mathrm{i}\eta$, with
\begin{eqnarray*}
\Phi_{i}(z)=\hat{\Phi}_{i}\big{(}\phi(z)\big{)},
\end{eqnarray*}
and $\phi(z)=\zeta$ is the inverse conformal mapping of $z=\psi(\zeta)$ which conformally maps $S$ to the interior region of a unit circle $|\zeta|<1$ and $\partial S$ to its circumference $|\zeta|=1$.\\

\noindent Let
\begin{eqnarray*}
\hat{\lambda}_{1}(z)=-B_{b}(z)-\mathrm{i}A_{b}(z),z\in\partial S, \\
\hat{\lambda}_{2}(z)=-D_{b}(z)+\mathrm{i}C_{b}(z),z\in\partial S.
\end{eqnarray*}
Denote
$$
\tilde{g}_{i}(\zeta)=g_{i}[\psi(\zeta)],
\tilde{\lambda}_{i}(\zeta)=\hat{\lambda}_{i}[\psi(\zeta)], \tilde{\tau}_{i}(\zeta)=\tau_{i}[\psi(\zeta)],
$$
and
\begin{eqnarray*}
& \tilde{q}_{i}(\zeta)=q_{i}[\psi(\zeta)]=\arg \left[\tilde{\lambda}_{i}(\zeta) \mathrm{e}^{\tilde{\tau}_{i}(\zeta)}\right]+m_{i} \arg [\psi(\zeta)],  \\
& \tilde{\chi}_{i}(\zeta)=\frac{1}{2 \pi} \int_{|t|=1} \tilde{q}_{i}(t) \frac{t+\zeta}{t-\zeta} \frac{\mathrm{d} t}{t},
\end{eqnarray*}
the integers $m_{i}$ are so chosen that every branch of $\tilde{q}_{i}(\zeta)$ is a single valued function on $\{\zeta:|\zeta|=1\}$.\\

\noindent If $m_{i} \geq 0, \tilde{\Phi}_{i}(\zeta)$ are given by
\begin{eqnarray*}
\tilde{\Phi}_{i}(\zeta)=\frac{\zeta^{m_{i}} \mathrm{e}^{-\tilde{\chi}_{i}(\zeta)}}{2 \pi \mathrm{i}} \int_{|t|=1} \hat{g}_{i}(t) \frac{t+\zeta}{t-\zeta} \frac{\mathrm{d} t}{t}+\mathrm{e}^{-\tilde{\chi}_{i}(\zeta)} \sum_{k=0}^{2 m_{i}} c_{k} \zeta^{k},
\end{eqnarray*}
where
\begin{eqnarray*}
\hat{g}_{i}(\zeta)=\frac{\tilde{g}_{i}(\zeta) \mathrm{e}^{\tilde{p}_{i}(\zeta)}}{\left|\tilde{\lambda}_{i}(\zeta) \mathrm{e}^{\tilde{\tau}_{i}(\zeta)}\right|},
\end{eqnarray*}
with $\tilde{p}_{i}(\zeta)$ being the real part of $\tilde{\chi}_{i}(\zeta)$. $c_{0}, c_{1}, \ldots, c_{2 m_{i}}$ are constants satisfying
\begin{eqnarray*}
c_{2 m_{i-k}}=-\bar{c}_{k}, \quad k=0,1, \ldots, m_{i}.
\end{eqnarray*}

\noindent If $m_{i}<0, \tilde{\Phi}_{i}(\zeta)$ are expressed by
\begin{eqnarray*}
\tilde{\Phi}_{i}(\zeta)=\frac{\mathrm{e}^{-\tilde{\chi}_{i}(\zeta)}}{\pi \mathrm{i}} \int_{|t|=1} \frac{\hat{g}_{i}(t) t^{m_{i}} \mathrm{d} t}{t-\zeta},
\end{eqnarray*}
when and only when
\begin{eqnarray*}
\int_{0}^{2 \pi} \hat{g}_{i}\left(e^{\mathrm{i} \theta}\right) \mathrm{e}^{-k \mathrm{i} \theta} \mathrm{d} \theta=0, \quad k=0, \ldots,-m_{i}+1.
\end{eqnarray*}
\end{lem}

\begin{Cor}\label{Cor4.2}
The RHBVP for the bi-axial meta-monogenic functions in a general case, a bi-axial, monogenic, continuous differentiable, $\mathbb{R}_{0,n}$-valued function $\varphi$, satisfying
\begin{eqnarray}\label{3.15}
\left\{\begin{array}{ll} (\mathcal{D}-\alpha)\varphi_(\underline{x})=0,&\ \ \underline{x}\in \Omega,\\
\mathrm{Re}\{\lambda_1(t)\varphi_2(t)\}=g_1(t),&\ \ t\in\partial \Omega,\\
\mathrm{Re}\{\lambda_2(t)\varphi_2(t)\}=g_2(t),&\ \ t\in\partial \Omega,
\end{array}\right.
\end{eqnarray}
is solvable.\\

\noindent Taking $\beta=\sum_{j=1}^{n}\alpha_jx_j$, we obtain
\begin{eqnarray*}
(\mathcal{D}-\alpha)\varphi=0{\rm\ \ \Leftrightarrow\ }\mathcal{D}(\mathrm{e}^{-\beta}\varphi_1)+\mathcal{D}(\mathrm{e}^{-\beta}\varphi_2)=0.
\end{eqnarray*}
Then, the problem \eqref{3.15} is equivalent to
\begin{eqnarray*}
\left\{\begin{array}{ll} \mathcal{D}(\mathrm{e}^{-\beta}\varphi_1)=0,&\quad \underline{x}\in \Omega,\\
\mathcal{D}(\mathrm{e}^{-\beta}\varphi_2)=0,&\quad \underline{x}\in \Omega,\\
\mathrm{Re}\{\lambda_1(t)\mathrm{e}^{-\beta}\varphi_1\}=\mathrm{e}^{-\beta}g_1(t),&\quad t\in\partial \Omega,\\
\mathrm{Re}\{\lambda_2(t)\mathrm{e}^{-\beta}\varphi_2\}=\mathrm{e}^{-\beta}g_2(t),&\quad t\in\partial \Omega.
\end{array}\right.
\end{eqnarray*}

\noindent Therefore the solution to the RHBVP \eqref{3.15} can be written as
\begin{eqnarray*}
\varphi(\underline{x})=e^{\beta}\big{(} \mathrm{Im}f_1(r,\rho)+\underline{\omega}\underline{\nu}\mathrm{Re}f_1(r,\rho)
+\underline{\omega}\mathrm{Im}f_2(r,\rho)+\underline{\nu}\mathrm{Re}f_2(r,\rho) \big{)},\ \underline{x}\in\Omega,
\end{eqnarray*}
where the function $f_i(i=1,2)$ is given by \eqref{3.11} in corresponding solvable conditions.
\end{Cor}

\medskip
\section{Riemann-Hilbert problems}

This section rigorously addresses Riemann-Hilbert boundary value problems (RHBVPs) for bi-axially symmetric poly-monogenic functions, the null solutions to iterated Dirac operators, defined over Clifford algebra-valued variable coefficient systems in $\mathbb{R}^{n}$. We first resolve the foundational case of these problems, including the Schwarz problem (a special RHBVP with constant coefficient 1), by leveraging the derived Almansi-type decomposition theorem to construct explicit integral representations tailored to biaxial symmetry. Building on this framework, we then extend our analysis to the RHBVPs governing bi-axially symmetric meta-monogenic functions, defined as the null solutions to iterated perturbed Dirac operators with a vector wave number.\\

\noindent Let $\Omega\subset\mathbb{R}^n$ be a star-like sub-domain with the center at $0$ and be bi-axially symmetric. The RHBVP with Clifford-valued variable coefficients for poly-monogenic functions of biaxial type is: to find a function
$\varphi\in\mathcal{C}^{k}\big(\Omega,\mathbb{R}_{0,n}\big)\big(k\geq2,k\in\mathbb{N}\big)$, satisfying
\begin{equation}\label{4.1}
\left\{\begin{array}{ll}\mathcal{D}^k\varphi(\underline{x})=0,\ \ \underline{x}\in\Omega,\\
 \mathrm{Re}\{\lambda_1(\underline{t})(\mathcal{D}^{l}\varphi)(\underline{t})\}= g_{1,l}(\underline{t}),\ \ \underline{t}\in\partial\Omega,0\leq l\leq k-1,\  l\in\mathbb{N},\\
 \mathrm{Re}\{\lambda_2(\underline{t})(\mathcal{D}^{l}\varphi)(\underline{t})\}= g_{2,l}(\underline{t}),\ \ \underline{t}\in\partial\Omega,0\leq l\leq k-1,\  l\in\mathbb{N},
\end{array}\right.
\end{equation}
where $\lambda_1(\underline{t})=A_b(r,\rho)+\underline{\omega}\underline{\nu}B_b(r,\rho)$ is of the first bi-axial type, $\lambda_2(\underline{t})
=\underline{\omega}C_b(r,\rho)+\underline{\nu}D_b(r,\rho)$ is of the second bi-axial type, and $\lambda_{i}\in\mathcal{H}^{\mu}(\partial\Omega,\mathbb{R}_{0,n}),g_{i}\in\mathcal{H}^{\mu}(\partial \Omega,\mathbb{R})(i=1,2)$.

\begin{Th}\label{Th4.1}
The RHBVP \eqref{4.1} is solvable, and its unique solution is given by
$$
\varphi(\underline{x})=\sum_{j=0}^{k-1} \underline{x}^{j} \varphi_{j}(\underline{x}), \ \underline{x} \in \Omega.
$$

\noindent For $j=0,1,\ldots,k-1$, $\varphi_j$ are expressed as
$$
\left\{\begin{array}{ll}
\varphi_{0}=\tilde{\varphi}_{0}-\sum_{j=1}^{k-1}\underline{x}^{j}\varphi_{j}, & j=0, \\
\varphi_{1}=C_{1, 1}^{-1} I_{\frac{n+1}{2}}
\left(\tilde{\varphi}_{1}-\sum_{j=2}^{k-1}\mathfrak{C}_{1, j} \underline{x}^{j-1}E_{\frac{n+1}{2}+\left[\frac{j-1}{2}\right]} \ldots E_{\frac{n+1}{2}+\left[\frac{j}{2}\right]-1} \varphi_{j} \right), & j=1, \\
\varphi_{l}=C_{l, l}^{-1} I_{\frac{n+1}{2}+\left[\frac{l}{2}\right]-1} \ldots I_{\frac{n+1}{2}}
\left(\tilde{\varphi}_{l}-\sum_{j=l+1}^{k-1}\mathfrak{C}_{l, j} \underline{x}^{j-l}E_{\frac{n+1}{2}+\left[\frac{j-l}{2}\right]} \ldots E_{\frac{n+1}{2}+\left[\frac{j}{2}\right]-1} \varphi_{j} \right), & 2 \leq l \leq k-1,
\end{array}\right.
$$
where each $\tilde{\varphi}_{l}(l=0,1,\ldots,k-1)$ is written by
$$
\tilde{\varphi}_{l}(\underline{x})=\mathrm{Im}\left\{f_{1,l}(r,\rho)\right\}+\underline{\omega\nu}\mathrm{Re}\left\{f_{1,l}(r,\rho)\right\}
+\underline{\omega}\mathrm{Im}\left\{f_{2,l}(r,\rho)\right\}+\underline{\nu}\mathrm{Re}\left\{f_{2,l}(r,\rho)\right\},
$$
where $f_{1,l},f_{2,l}$ are given by expressions \eqref{3.11},\eqref{3.12},\eqref{3.13} in Lemma \ref{3.4} and the corresponding functions $g_{1,l},g_{2,l}$.
\end{Th}

\noindent\p\ Since $\mathcal{D}^{k} \varphi=0$, by applying Corollary \ref{Cor3.1}, we know that there exist unique functions $\varphi_{j}$ satisfying $\mathcal{D} \varphi_{j}=0(j=0,1,2, \ldots, k-1)$, such that
$$
\varphi=\sum_{j=0}^{k-1} \underline{x}^{j} \varphi_{j}, \underline{x} \in \Omega,
$$
where each $\varphi_j$ is given by \eqref{3.7}.
By Lemma \ref{lem3.3}, for $l \in \mathbb{N}$ and $l \leq j$, we obtain
$$
\mathcal{D}^{l} \varphi=\sum_{j=0}^{k-1} \mathcal{D}^{l}\left(\underline{x}^{j} \varphi_{j}\right)=\sum_{j=l}^{k-1} \mathfrak{C}_{l, j} \underline{x}^{j-l} E_{\frac{n+1}{2}+\left[\frac{j-l}{2}\right]} \ldots E_{\frac{n+1}{2}+\left[\frac{j}{2}\right]-1} \varphi_{j}.
$$

\noindent Now, the RHBVP \eqref{4.1} is equivalent to the problem
\begin{eqnarray}\label{4.2}
\left\{\begin{array}{ll}
\mathcal{D} \varphi_{j}(\underline{x})=0,\ \ j=0,1,2, \ldots, k-1,\\
\mathrm{Re}\left\{\lambda_1(\underline{t})\sum_{j=0}^{k-1} \underline{t}^{j}\varphi_{j}(\underline{t}) \right\} = g_{1,0}(\underline{t}),\\
\mathrm{Re}\left\{\lambda_2(\underline{t})\sum_{j=0}^{k-1} \underline{t}^{j}\varphi_{j}(\underline{t}) \right\} = g_{2,0}(\underline{t}),\\
\mathrm{Re}\left\{\lambda_1(\underline{t})\sum_{j=1}^{k-1} \mathfrak{C}_{1, j} \underline{t}^{j-1}E_{\frac{n+1}{2}+\left[\frac{j-1}{2}\right]} \ldots E_{\frac{n+1}{2}+\left[\frac{j}{2}\right]-1} \varphi_{j}(\underline{t}) \right\} = g_{1,1}(\underline{t}),\\
\mathrm{Re}\left\{\lambda_2(\underline{t})\sum_{j=1}^{k-1} \mathfrak{C}_{1, j} \underline{t}^{j-1}E_{\frac{n+1}{2}+\left[\frac{j-1}{2}\right]} \ldots E_{\frac{n+1}{2}+\left[\frac{j}{2}\right]-1} \varphi_{j}(\underline{t}) \right\} = g_{2,1}(\underline{t}),\\
\vdots \\
\mathrm{Re}\left\{\lambda_1(\underline{t}) \mathfrak{C}_{k-1, k-1} E_{\frac{n+1}{2}} \ldots E_{\frac{n+1}{2}+\left[\frac{k-1}{2}\right]-1} \varphi_{k-1}(\underline{t}) \right\} = g_{1,k-1}(\underline{t}),\\
\mathrm{Re}\left\{\lambda_2(\underline{t}) \mathfrak{C}_{k-1, k-1} E_{\frac{n+1}{2}} \ldots E_{\frac{n+1}{2}+\left[\frac{k-1}{2}\right]-1} \varphi_{k-1}(\underline{t}) \right\} = g_{2,k-1}(\underline{t}).
\end{array}\right.
\end{eqnarray}

\noindent Next, we consider the RHBVP for the case $k-1$:
$$
\left\{\begin{array}{l}
\mathcal{D} {\varphi}_{k-1}(\underline{x})=0, \\
\mathrm{Re}\left\{\lambda_1(\underline{t}) \mathfrak{C}_{k-1, k-1} E_{\frac{n+1}{2}} \ldots E_{\frac{n+1}{2}+\left[\frac{k-1}{2}\right]-1}\varphi_{k-1}(\underline{t}) \right\} = g_{1,k-1}(\underline{t}),\\
\mathrm{Re}\left\{\lambda_2(\underline{t}) \mathfrak{C}_{k-1, k-1} E_{\frac{n+1}{2}} \ldots E_{\frac{n+1}{2}+\left[\frac{k-1}{2}\right]-1}\varphi_{k-1}(\underline{t}) \right\} = g_{2,k-1}(\underline{t}). \\
\end{array}\right.
$$

\noindent Denote
$$
\tilde{\varphi}_{k-1}(\underline{x})= \mathfrak{C}_{k-1, k-1} E_{\frac{n+1}{2}} \ldots E_{\frac{n+1}{2}+\left[\frac{k-1}{2}\right]-1}\varphi_{k-1}(\underline{x}),
$$
then $\tilde{\varphi}_{k-1}$ is a monogenic function of biaxial type. Using Lemma \ref{3.4}, we can obtain that
$$
\tilde{\varphi}_{k-1}(\underline{x})=\mathrm{Im}\left\{f_{1,k-1}(r,\rho)\right\}+\underline{\omega\nu}\mathrm{Re}\left\{f_{1,k-1}(r,\rho)\right\}
+\underline{\omega}\mathrm{Im}\left\{f_{2,k-1}(r,\rho)\right\}+\underline{\nu}\mathrm{Re}\left\{f_{2,k-1}(r,\rho)\right\},
$$
where $f_{1,k-1},f_{2,k-1}$ are given by \eqref{3.11}. Therefore, applying Lemma \ref{lem3.1}, $\varphi_{k-1}$ is given by
$$
\varphi_{k-1}=\mathfrak{C}_{k-1, k-1}^{-1} I_{\frac{n+1}{2}+\left[\frac{k-1}{2}\right]-1}\ldots I_{\frac{n+1}{2}} \tilde{\varphi}_{k-1}.
$$

\noindent Then, we can proceed in an inductive way. Applying Lemmas \ref{3.1}, \ref{3.4}, we consider the RHBVP for case $k-2$:
$$
\left\{\begin{array}{l}
\mathcal{D} {\varphi}_{k-2}(\underline{x})=0, \\
\mathrm{Re}\left\{\lambda_1 \left(\mathfrak{C}_{k-2, k-2} E_{\frac{n+1}{2}} \ldots E_{\frac{n+1}{2}+\left[\frac{k-2}{2}\right]-1}\varphi_{k-2}
+\mathfrak{C}_{k-2, k-1}\underline{t} E_{\frac{n+1}{2}+\frac{1}{2}} \ldots E_{\frac{n+1}{2}+\left[\frac{k-1}{2}\right]-1}\varphi_{k-1} \right) \right\} = g_{1,k-2},\\
\mathrm{Re}\left\{\lambda_2 \left(\mathfrak{C}_{k-2, k-2} E_{\frac{n+1}{2}} \ldots E_{\frac{n+1}{2}+\left[\frac{k-2}{2}\right]-1}\varphi_{k-2}
+\mathfrak{C}_{k-2, k-1}\underline{t} E_{\frac{n+1}{2}+\frac{1}{2}} \ldots E_{\frac{n+1}{2}+\left[\frac{k-1}{2}\right]-1}\varphi_{k-1} \right) \right\} = g_{2,k-2}.
\end{array}\right.
$$

\noindent Denote
$$
\tilde{\varphi}_{k-2}= \mathfrak{C}_{k-2, k-2} E_{\frac{n+1}{2}} \ldots E_{\frac{n+1}{2}+\left[\frac{k-2}{2}\right]-1}\varphi_{k-2}
+\mathfrak{C}_{k-2, k-1}\underline{x} E_{\frac{n+1}{2}+\frac{1}{2}} \ldots E_{\frac{n+1}{2}+\left[\frac{k-1}{2}\right]-1}\varphi_{k-1},
$$
then, $\tilde{\varphi}_{k-2}$ is also a monogenic function of biaxial type.

\noindent In a similar way, we have
$$
\tilde{\varphi}_{k-2}(\underline{x})=\mathrm{Im}\left\{f_{1,k-2}(r,\rho)\right\}+\underline{\omega\nu}\mathrm{Re}\left\{f_{1,k-2}(r,\rho)\right\}
+\underline{\omega}\mathrm{Im}\left\{f_{2,k-2}(r,\rho)\right\}+\underline{\nu}\mathrm{Re}\left\{f_{2,k-2}(r,\rho)\right\}.
$$

\noindent As above we prove that RHBVP for case $k-2$ has the solution
$$
\varphi_{k-2}=\mathfrak{C}_{k-2, k-2}^{-1} I_{\frac{n+1}{2}+\left[\frac{k-2}{2}\right]-1}\ldots I_{\frac{n+1}{2}}
\left( \tilde{\varphi}_{k-1}-\mathfrak{C}_{k-2, k-1}\underline{t} E_{\frac{n+1}{2}+\frac{1}{2}} \ldots E_{\frac{n+1}{2}+\left[\frac{k-1}{2}\right]-1}\varphi_{k-1}\right).
$$

\noindent By induction for $2 \leq l \leq k-1$, the boundary value problem
$$
\left\{\begin{array}{ll}
\mathcal{D} \varphi_{l}(\underline{x})=0, \\
\mathrm{Re}\left\{\lambda_1(\underline{t})\sum_{j=l}^{k-1} \mathfrak{C}_{l, j} \underline{t}^{j-l}E_{\frac{n+1}{2}+\left[\frac{j-l}{2}\right]} \ldots E_{\frac{n+1}{2}+\left[\frac{j}{2}\right]-1} \varphi_{j}(\underline{t}) \right\} = g_{1,l}(\underline{t}),\\
\mathrm{Re}\left\{\lambda_2(\underline{t})\sum_{j=l}^{k-1} \mathfrak{C}_{l, j} \underline{t}^{j-l}E_{\frac{n+1}{2}+\left[\frac{j-l}{2}\right]} \ldots E_{\frac{n+1}{2}+\left[\frac{j}{2}\right]-1} \varphi_{j}(\underline{t}) \right\} = g_{2,l}(\underline{t})
\end{array}\right.
$$
has the solution
$$
\tilde{\varphi}_{l}(\underline{x})=\mathrm{Im}\left\{f_{1,l}(r,\rho)\right\}+\underline{\omega\nu}\mathrm{Re}\left\{f_{1,l}(r,\rho)\right\}
+\underline{\omega}\mathrm{Im}\left\{f_{2,l}(r,\rho)\right\}+\underline{\nu}\mathrm{Re}\left\{f_{2,l}(r,\rho)\right\},
$$
with
$$
\tilde{\varphi}_{l}=\sum_{j=l}^{k-1} \mathfrak{C}_{l, j} \underline{x}^{j-l}E_{\frac{n+1}{2}+\left[\frac{j-l}{2}\right]} \ldots E_{\frac{n+1}{2}+\left[\frac{j}{2}\right]-1} \varphi_{j}.
$$

\noindent Therefore, $\varphi_l$ is expressed by
$$
\varphi_{l}=C_{l, l}^{-1} I_{\frac{n+1}{2}+\left[\frac{l}{2}\right]-1} \ldots I_{\frac{n+1}{2}}
\left(\tilde{\varphi}_{l}-\sum_{j=l+1}^{k-1}\mathfrak{C}_{l, j} \underline{x}^{j-l}E_{\frac{n+1}{2}+\left[\frac{j-l}{2}\right]} \ldots E_{\frac{n+1}{2}+\left[\frac{j}{2}\right]-1} \varphi_{j} \right).
$$

\noindent For $l=1$, the boundary value problem
$$
\left\{\begin{array}{ll}
\mathcal{D} \varphi_{l}(\underline{x})=0, \\
\mathrm{Re}\left\{\lambda_1(\underline{t})\sum_{j=1}^{k-1} \mathfrak{C}_{1, j} \underline{t}^{j-1}E_{\frac{n+1}{2}+\left[\frac{j-1}{2}\right]} \ldots E_{\frac{n+1}{2}+\left[\frac{j}{2}\right]-1} \varphi_{j}(\underline{t}) \right\} = g_{1,1}(\underline{t}),\\
\mathrm{Re}\left\{\lambda_2(\underline{t})\sum_{j=1}^{k-1} \mathfrak{C}_{1, j} \underline{t}^{j-1}E_{\frac{n+1}{2}+\left[\frac{j-1}{2}\right]} \ldots E_{\frac{n+1}{2}+\left[\frac{j}{2}\right]-1} \varphi_{j}(\underline{t}) \right\} = g_{2,1}(\underline{t})
\end{array}\right.
$$
has the solution
$$
\tilde{\varphi}_{1}(\underline{x})=\mathrm{Im}\left\{f_{1,1}(r,\rho)\right\}+\underline{\omega\nu}\mathrm{Re}\left\{f_{1,1}(r,\rho)\right\}
+\underline{\omega}\mathrm{Im}\left\{f_{2,1}(r,\rho)\right\}+\underline{\nu}\mathrm{Re}\left\{f_{2,1}(r,\rho)\right\},
$$
with
$$
\tilde{\varphi}_{1}=\sum_{j=1}^{k-1} \mathfrak{C}_{1, j} \underline{x}^{j-1}E_{\frac{n+1}{2}+\left[\frac{j-1}{2}\right]} \ldots E_{\frac{n+1}{2}+\left[\frac{j}{2}\right]-1} \varphi_{j}.
$$

\noindent Therefore, $\varphi_1$ is expressed by
$$
\varphi_{1}=C_{1, 1}^{-1} I_{\frac{n+1}{2}}
\left(\tilde{\varphi}_{1}-\sum_{j=2}^{k-1}\mathfrak{C}_{1, j} \underline{x}^{j-1}E_{\frac{n+1}{2}+\left[\frac{j-1}{2}\right]} \ldots E_{\frac{n+1}{2}+\left[\frac{j}{2}\right]-1} \varphi_{j} \right).
$$

\noindent Finally, for $l=0$, in the same way we can show that the boundary value problem
$$
\left\{\begin{array}{ll}
\mathcal{D} \varphi_{0}(\underline{x})=0, \\
\mathrm{Re}\left\{\lambda_1(\underline{t})\sum_{j=0}^{k-1} \underline{t}^{j}\varphi_{j}(\underline{t}) \right\} = g_{1,0}(\underline{t}),\\
\mathrm{Re}\left\{\lambda_2(\underline{t})\sum_{j=0}^{k-1} \underline{t}^{j}\varphi_{j}(\underline{t}) \right\} = g_{2,0}(\underline{t})
\end{array}\right.
$$
has the solution
$$
\tilde{\varphi}_{0}(\underline{x})=\mathrm{Im}\left\{f_{1,0}(r,\rho)\right\}+\underline{\omega\nu}\mathrm{Re}\left\{f_{1,0}(r,\rho)\right\}
+\underline{\omega}\mathrm{Im}\left\{f_{2,0}(r,\rho)\right\}+\underline{\nu}\mathrm{Re}\left\{f_{2,0}(r,\rho)\right\},
$$
with
$$
\tilde{\varphi}_{0}=\sum_{j=0}^{k-1} \underline{x}^{j}\varphi_{j}.
$$

\noindent Therefore, $\varphi_0$ is expressed by
$$
\varphi_{0}=\tilde{\varphi}_{0}-\sum_{j=1}^{k-1}\underline{x}^{j}\varphi_{j}.\\
$$

\noindent Combining the above terms, the RHBVP \eqref{4.1} is solvable and its unique solution is given by
$$
\varphi(\underline{x})=\sum_{j=0}^{k-1} \underline{x}^{j} \varphi_{j}(\underline{x}), \ \underline{x} \in \Omega,
$$
where all $\varphi_j(j=0,1,\ldots,k-1)$ are given as above.\\

\noindent The proof is completed.
$\hfill\square$\\

\noindent In what follows let us continue to consider RHBVP forthe null solution to the perturbed iterated Dirac operator with a vector wave number.\\

\begin{Th}
Let $\Omega\subset\mathbb{R}^n$ be a star-like sub-domain with the center at $0$ and be bi-axially symmetric. The RHBVP with Clifford-valued variable coefficients for function which is the null solution to the perturbed iterated Dirac operator with a vector wave number and of bi-axial type, is to find a function
$\varphi\in\mathcal{C}^{k}\big(\Omega,\mathbb{R}_{0,n}\big)\big(k\geq2,k\in\mathbb{N}\big)$, satisfying
\begin{equation}\label{4.3}
\left\{\begin{array}{ll}(\mathcal{D}-\alpha)^k\varphi(\underline{x})=0,\ \ \underline{x}\in\Omega,\\
 \mathrm{Re}\{\lambda_1(\underline{t})(\mathcal{D}-\alpha)^{l}\varphi (\underline{t})\}= g_{1,l}(\underline{t}),\ \ \underline{t}\in\partial\Omega,0\leq l\leq k-1,\  l\in\mathbb{N},\\
 \mathrm{Re}\{\lambda_2(\underline{t})(\mathcal{D}-\alpha)^{l}\varphi(\underline{t})\}= g_{2,l}(\underline{t}),\ \ \underline{t}\in\partial\Omega,0\leq l\leq k-1,\  l\in\mathbb{N},
\end{array}\right.
\end{equation}
where $\alpha=\sum_{j=1}^{n}e_{j}\alpha_{j},\alpha_{j}\in\mathbb{R}$, $\lambda_1(\underline{t})=A_b(r,\rho)+\underline{\omega}\underline{\nu}B_b(r,\rho)$ is of the first biaxial type, $\lambda_2(\underline{t})
=\underline{\omega}C_b(r,\rho)+\underline{\nu}D_b(r,\rho)$ is of the second biaxial type,  and $\lambda_{i}\in\mathcal{H}^{\mu}(\partial\Omega,\mathbb{R}_{0,n}),g_{i}\in\mathcal{H}^{\mu}(\partial \Omega,\mathbb{R})(i=1,2)$.

\noindent The RHBVP is solvable and its unique solution is given by
$$
\varphi(\underline{x})=\mathrm{e}^{\beta} \sum_{j=0}^{k-1} \underline{x}^{j} \varphi_{j}(\underline{x}), \ \underline{x} \in \Omega,
$$
with $\beta=\sum_{j=1}^{n}\alpha_{j}x_{j}$.
For $j=0,1,\ldots,k-1$, $\varphi_j$ are expressed as
$$
\left\{\begin{array}{ll}
\varphi_{0}=\tilde{\varphi}_{0}-\mathrm{e}^{\beta} \sum_{j=1}^{k-1}\underline{x}^{j}\varphi_{j}, & j=0, \\
\varphi_{1}=C_{1, 1}^{-1} \mathrm{e}^{-\beta} I_{\frac{n+1}{2}}
\left(\tilde{\varphi}_{1}-\mathrm{e}^{\beta} \sum_{j=2}^{k-1}\mathfrak{C}_{1, j} \underline{x}^{j-1}E_{\frac{n+1}{2}+\left[\frac{j-1}{2}\right]} \ldots E_{\frac{n+1}{2}+\left[\frac{j}{2}\right]-1} \varphi_{j} \right), & j=1, \\
\varphi_{l}=C_{l, l}^{-1} \mathrm{e}^{-\beta} I_{\frac{n+1}{2}+\left[\frac{l}{2}\right]-1} \ldots I_{\frac{n+1}{2}}
\left(\tilde{\varphi}_{l}-\mathrm{e}^{\beta} \sum_{j=l+1}^{k-1}\mathfrak{C}_{l, j} \underline{x}^{j-l}E_{\frac{n+1}{2}+\left[\frac{j-l}{2}\right]} \ldots E_{\frac{n+1}{2}+\left[\frac{j}{2}\right]-1} \varphi_{j} \right), & 2 \leq l \leq k-1,
\end{array}\right.
$$
where each $\tilde{\varphi}_{l}(l=0,1,\ldots,k-1)$ is written by
$$
\tilde{\varphi}_{l}(\underline{x})=\mathrm{Im}\left\{f_{1,l}(r,\rho)\right\}+\underline{\omega\nu}\mathrm{Re}\left\{f_{1,l}(r,\rho)\right\}
+\underline{\omega}\mathrm{Im}\left\{f_{2,l}(r,\rho)\right\}+\underline{\nu}\mathrm{Re}\left\{f_{2,l}(r,\rho)\right\},
$$
where $f_{1,l},f_{2,l}$ are given by expressions \eqref{3.11},\eqref{3.12},\eqref{3.13} in Lemma \ref{3.4} and the corresponding functions $g_{1,l},g_{2,l}$.
\end{Th}

\noindent\p\ Since $(\mathcal{D}-\alpha)^{k} \varphi=0$, by applying Lemma \ref{lem3.2} we know that there exist unique functions $\varphi_{j}$ satisfying $\mathcal{D} \varphi_{j}=0(j=0,1,2, \ldots, k-1)$, such that
$$
\varphi=\mathrm{e}^{\beta}\sum_{j=0}^{k-1} \underline{x}^{j} \varphi_{j}, \underline{x} \in \Omega,
$$
where $\beta=\sum_{j=1}^{n}\alpha_{j}x_{j}$ and each $\varphi_j$ is given by \eqref{8}.
For $l \in \mathbb{N}$ and $l \leq j$, we obtain
$$
(\mathcal{D}-\alpha)^{l} \varphi=\sum_{j=0}^{k-1} (\mathcal{D}-\alpha)^{l}\left( \mathrm{e}^{\beta} \underline{x}^{j} \varphi_{j}\right)=\mathrm{e}^{\beta} \sum_{j=l}^{k-1} \mathfrak{C}_{l, j} \underline{x}^{j-l} E_{\frac{n+1}{2}+\left[\frac{j-l}{2}\right]} \ldots E_{\frac{n+1}{2}+\left[\frac{j}{2}\right]-1} \varphi_{j}.
$$

\noindent Therefore, the RHBVP \eqref{4.3} is equivalent to the problem
\begin{eqnarray}
\left\{\begin{array}{ll}
\mathcal{D} \varphi_{j}(\underline{x})=0,\ \ j=0,1,2, \ldots, k-1,\\
\mathrm{Re}\left\{\lambda_1(\underline{t})\mathrm{e}^{\beta} \sum_{j=0}^{k-1} \underline{t}^{j}\varphi_{j}(\underline{t}) \right\} = g_{1,0}(\underline{t}),\\
\mathrm{Re}\left\{\lambda_2(\underline{t})\mathrm{e}^{\beta} \sum_{j=0}^{k-1} \underline{t}^{j}\varphi_{j}(\underline{t}) \right\} = g_{2,0}(\underline{t}),\\
\mathrm{Re}\left\{\lambda_1(\underline{t})\mathrm{e}^{\beta} \sum_{j=1}^{k-1} \mathfrak{C}_{1, j} \underline{t}^{j-1}E_{\frac{n+1}{2}+\left[\frac{j-1}{2}\right]} \ldots E_{\frac{n+1}{2}+\left[\frac{j}{2}\right]-1} \varphi_{j}(\underline{t}) \right\} = g_{1,1}(\underline{t}),\\
\mathrm{Re}\left\{\lambda_2(\underline{t})\mathrm{e}^{\beta} \sum_{j=1}^{k-1} \mathfrak{C}_{1, j} \underline{t}^{j-1}E_{\frac{n+1}{2}+\left[\frac{j-1}{2}\right]} \ldots E_{\frac{n+1}{2}+\left[\frac{j}{2}\right]-1} \varphi_{j}(\underline{t}) \right\} = g_{2,1}(\underline{t}),\\
\vdots \\
\mathrm{Re}\left\{\lambda_1(\underline{t})\mathrm{e}^{\beta} \mathfrak{C}_{k-1, k-1} E_{\frac{n+1}{2}} \ldots E_{\frac{n+1}{2}+\left[\frac{k-1}{2}\right]-1} \varphi_{k-1}(\underline{t}) \right\} = g_{1,k-1}(\underline{t}),\\
\mathrm{Re}\left\{\lambda_2(\underline{t})\mathrm{e}^{\beta} \mathfrak{C}_{k-1, k-1} E_{\frac{n+1}{2}} \ldots E_{\frac{n+1}{2}+\left[\frac{k-1}{2}\right]-1} \varphi_{k-1}(\underline{t}) \right\} = g_{2,k-1}(\underline{t}).
\end{array}\right.
\end{eqnarray}
Similar to the method in Theorem \ref{Th4.1}, we can obtain the solution
$$
\varphi(\underline{x})=\mathrm{e}^{\beta} \sum_{j=0}^{k-1} \underline{x}^{j} \varphi_{j}(\underline{x}), \ \underline{x} \in \Omega,
$$
where all $\varphi_j$ are expressed as
$$
\left\{\begin{array}{ll}
\varphi_{0}=\tilde{\varphi}_{0}-\mathrm{e}^{\beta} \sum_{j=1}^{k-1}\underline{x}^{j}\varphi_{j}, & j=0, \\
\varphi_{1}=C_{1, 1}^{-1} \mathrm{e}^{-\beta} I_{\frac{n+1}{2}}
\left(\tilde{\varphi}_{1}-\mathrm{e}^{\beta} \sum_{j=2}^{k-1}\mathfrak{C}_{1, j} \underline{x}^{j-1}E_{\frac{n+1}{2}+\left[\frac{j-1}{2}\right]} \ldots E_{\frac{n+1}{2}+\left[\frac{j}{2}\right]-1} \varphi_{j} \right), & j=1, \\
\varphi_{l}=C_{l, l}^{-1} \mathrm{e}^{-\beta} I_{\frac{n+1}{2}+\left[\frac{l}{2}\right]-1} \ldots I_{\frac{n+1}{2}}
\left(\tilde{\varphi}_{l}-\mathrm{e}^{\beta} \sum_{j=l+1}^{k-1}\mathfrak{C}_{l, j} \underline{x}^{j-l}E_{\frac{n+1}{2}+\left[\frac{j-l}{2}\right]} \ldots E_{\frac{n+1}{2}+\left[\frac{j}{2}\right]-1} \varphi_{j} \right), & 2 \leq l \leq k-1.
\end{array}\right.
$$
$\hfill\square$

\begin{Rem}
Under H\"{o}lder continuous boundary data conditions, we investigate Riemann-Hilbert boundary value problems (RHBVPs) for bi-axially symmetric poly-monogenic functions, the null solutions to the iterated Dirac equation, governed by Clifford-algebra-valued variable coefficients. For the associated Schwarz problems, characterized by constant coefficients $\lambda_1=1$ and $\lambda_2=1$, our methodology retains validity, demonstrating consistency with classical results. Notably, the case $k=1$ in RHBVP \eqref{3.9} emerges as a specialization of the generalized formulation RHBVP \eqref{4.1}, establishing a hierarchical framework for analyzing iterative operator hierarchies. Furthermore, when boundary data resides in
$\mathcal{L}_p(p<1<+\infty)$, analogous existence and uniqueness properties extend to bi-axially symmetric null solutions to $\left(\mathcal{D} - \alpha\right)^k\varphi=0, k \geq2, \alpha=\sum\limits_{j=1}^{n}e_{j}\alpha_{j},\alpha_{j}\in\mathbb{R}$, underscoring the robustness of our approach across function space topologies. This unified treatment bridges low-order and higher-order perturbed Dirac systems while preserving geometric symmetry constraints, thereby circumventing redundancies in traditional case-by-case analyses.

\end{Rem}

\medskip
\section*{Competing interests}
The authors declare that they have no conflict of interest.

\section*{Data Availability}
Data sharing not applicable to this article as no datasets were generated or analyzed during
the current study.

%\section*{Funding}
%This work was supported by National Natural Science
%Foundation of China (11601525).

%\section*{Authors' contributions}
%Fuli He and Yajun Hu wrote the main manuscript text. All authors reviewed the manuscript.

%\section*{Acknowledgement}

%The author would like  to acknowledge his sincere thanks to the %Editor and the reviewers  for many valuable comments and %suggestions.

\end{document}